\documentclass[12pt,reqno]{amsart}
\usepackage{amsmath,amssymb, graphics, amscd,latexsym }
\usepackage{epsfig}

\makeatletter

\newtheorem{Theorem}{Theorem}[section]
\newtheorem{Lemma}[Theorem]{Lemma}
\newtheorem{Corollary}[Theorem]{Corollary}
\newtheorem{Proposition}[Theorem]{Proposition}

\theoremstyle{remark}

\newtheorem{Example}[Theorem]{Example}

\begin{document}

\title{Framed and Oriented Links of Codimension $2$}

\author{Jianhua Wang}

\begin{abstract}  Sanderson [12] gave an isomorphism
    $\theta:\pi_m(\vee_{i=1}^r S^2_i)\longrightarrow
    \pi_m(\vee_{i=1}^{r+1}{\Bbb C}P_i^{\infty}).$
    In this paper we construct for any subset
    $\sigma\subset\{1,2,\cdots,r\}$ an isomorphism
    $\theta_{\sigma}$ from $\pi_m(\vee_{i=1}^r S^2_i)$ 
    to $\pi_m(\vee_{i=1}^{r+1}{\Bbb C}P_i^{\infty}).$
    The inclusion $S^2\vee S^2\hookrightarrow
    {\Bbb C}P^{\infty}\vee {\Bbb C}P^{\infty}$
    induces a homomorphism
    $f:\pi_m(S^2\vee S^2)\longrightarrow
    \pi_m({\Bbb C}P^{\infty}\vee {\Bbb C}P^{\infty})$.
    We also compute $f$ by evaluating $f$ on each factor
    in the Hilton splitting of $\pi_m (S^2\vee S^2)$, the 
    results in [12] concerning the case $m=4$ are generalized.
\vskip 4pt
\leftline{{\bf Keywords:} Framed and oriented links, 
          Seifert surface, Hilton splitting.}
\leftline{{\bf 2000 Mathematics Subject Classification:}
     55Q20, 57Q45}
\end{abstract}

\maketitle

\section{introduction}

         A link $M_1\sqcup M_2\sqcup\cdots\sqcup M_r\subset\Bbb R^m$
    is an ordered disjoint union of closed smooth submanifolds, $m\ge 3$.
    If these submanifolds are oriented then we call it an oriented link.
    A framing of a $k$-codimensional submanifold in $\Bbb R^m$ is a
    trivialization of its normal vector bundle, or equivalently, an ordered
    set of $k$ linearly independent normal vector fields. If every component
    of the link has a given framing then we call it a framed link.
    We will assume that the codimensions of the components are $2$,
    when nothing else is stated. The bordism groups of framed and
    oriented links with $r$ components of codimension $2$ in $\Bbb R^m$
    are denoted by $FL^2_{m,r}$ and
    $L^2_{m,r}$ respectively. The facts
\begin{eqnarray*}
              & FL^2_{m,r}\cong\pi_m(\vee_{i=1}^r S^2_i),
              & L^2_{m,r}\cong\pi_m(\vee_{i=1}^r {\Bbb C}P_i^{\infty})
\end{eqnarray*}    
    are well known by Pontryagin-Thom construction.

         The main purpose of this paper is to discuss the 
    relationship between
    $FL^2_{m,r}$ and $L^2_{m,r+1}$ and to compute the homomorphism
               $$f:\pi_m(S^2\vee S^2)\longrightarrow
                       \pi_m({\Bbb C}P^{\infty}\vee {\Bbb C}P^{\infty})
                       \cong\pi_m(S^2).$$
     Geometrically, $f$ is given by forgetting the framing of
     $[M_1\sqcup M_2]_{fr}\in FL^2_{m,2}$ and keeping the orientation
     determined by the framing, so we may call $f$ a forgetful
     homomorphism. 

         {\bf Main results and the organization} of this paper:
    In \S 2 we give a formula of reframing, suggested by Koschorke, 
    and discuss briefly the role of framing in the Hilton splitting. Let
    $\sigma\subset\{1,2,\cdots,r\}$ be any subset. We construct in \S 3 an
    isomorphism $\theta_{\sigma}:FL^2_{m,r}\longrightarrow
    L^2_{m,r+1}$. We recover Sanderson's isomorphism $\theta$ 
    by taking $\sigma=\{1,2,\cdots,r\}$. 
    The forgetful homomorphism $f$ is computed in \S 4
    by using the inverse of $\theta_{\sigma}$ and by choosing $\sigma=\phi$,
    in particular the following result of [12] is generalized: 
    in case $m=4$ it holds
    $f\circ\gamma_*\not=0$ for $\gamma=[\iota_1,\iota_2]$, 
    $[\iota_1,[\iota_1,\iota_2]]$ and $[\iota_2,[\iota_1,\iota_2]]$,
    $f\circ\gamma_*=0$ for $\gamma=\iota_1$ and $\iota_2$.     
   
         We work in the category of smooth manifolds.

         I'm grateful to my supervisor Prof.\ U. Koschorke for some ideas 
    and stimulating discussions, and to Prof. U. Kaiser for many 
    helps and useful suggestions. Thanks also to Prof.\ M.\ Heusener
    for nice talks. As the English version of this paper is finished
    in Bar-Ilan University I am also thankful to Prof.\ T.\ Nowik
    for being warmly hosted.

\section{framing}

          Let $M^{m-2}\subset \Bbb R^m$ be a closed submanifold with framing
     ${\mathcal F}=(v_1,v_2)$, and let
       $$ s:M\longrightarrow S^1\backsimeq SO(2) $$ 
     be a continuous map.   
     For $x\in M$ we can represent $s(x)\in S^1$
     by an orthonormal matrix
            $$\pmatrix a_{11}(x) & a_{12}(x) \\
                  a_{21}(x) & a_{22}(x)
                 \endpmatrix  . $$
     Define $v_1'$ and $v_2'$ by
           \begin{eqnarray}
            v_1'(x)&=&a_{11}(x)v_1(x)+a_{12}(x)v_2(x), \nonumber \\
            v_2'(x)&=&a_{21}(x)v_1(x)+a_{22}(x)v_2(x). \nonumber
          \end{eqnarray}
     $(v_1',v_2')$ is a new framing of $M$ and is denoted by $s{\mathcal F}$.
     Up to homotopy we may assume $s$ is differential.
     Let $-1\in S^1$ be a regular value of $s$ and consider
     $Z=s^{-1}(-1)\subset M$. Let $Z\times [-1,1]\subset M$ be a small tubular
     neighbourhood of $Z$ such that the positive direction of
     $[-1,1]$ is in agreement with the usual orientation
     of $S^1$. Up to homotopy $s{\mathcal F}$ is in fact 
     the $2\pi$-rotation of
     ${\mathcal F}$ in this neighbourhood, namely outside this
     neighbourhood it is the same as $(v_1,v_2)$ and inside it
\begin{eqnarray}
            v_1'(z,t)&=&v_1(z,t)\cos (t+1)\pi+v_2(z,t)\sin (t+1)\pi,   \\
            v_2'(z,t)&=&-v_1(z,t)\sin (t+1)\pi+v_2(z,t)\cos (t+1)\pi,
\end{eqnarray}
     where $z\in Z$ and $t\in [-1,1]$.

           Let $v_3$ be the normal vector field of $Z\subset M$, 
     determined by the
     orientation of $S^1$, provide $Z\subset\Bbb R^m$ with the framing
     $(v_1,v_2,v_3)$. We define 
         $$[Z,{\mathcal F},s]=[Z,(v_1,v_2,v_3)]\in\pi_m(S^3).$$
     Let $\eta:S^3\longrightarrow S^2$ be the Hopf map,
     Koschorke observed that $[Z,{\mathcal F},s]$ is the only obstruction to
     homotope ${\mathcal F}$ to $s{\mathcal F}$ and conjectured that
     $\eta_*[Z,{\mathcal F},s]$ and the difference
                $$[M,s{\mathcal F}]-[M,{\mathcal F}]\in\pi_m(S^2)$$ 
     should be related by some formula.

          It's well known that $M$ bounds a Seifert surface $F$.
     If the framing of $M$ induced by $F$ is homotopic to the 
     given framing $\mathcal F$ of
     $M$, then we say $M$ is $S$-framed.  In this case we have
     $[M,{\mathcal F}]=0$, and Turaev [15] proved
     $[M,s{\mathcal F}]=\eta_*[Z,{\mathcal F},s]$. Note that if $M$ is 
     oriented then the $S$-framing compatible with the orientation
     is unique up to homotopy.
\begin{Proposition} Let $M\subset\Bbb R^m$ be a closed submanifold
     with framing ${\mathcal F}=(v_1,v_2)$ and $s:M\longrightarrow S^1$
     be a map. Let $u{\mathcal F}=(-v_1,v_2)$ and $M^{sh}$ be a 
     small shift of $M$ along $v_1$ provided with the framing
     $s{\mathcal F}=(v'_1,v'_2)$. Then it holds
             $$[(M\sqcup M^{sh}),(u{\mathcal F}\sqcup s{\mathcal F})]=
                   \eta_*[Z,{\mathcal F},s].$$
\end{Proposition}
\begin{proof} Writing $(M\sqcup M^{sh})$ we mean it is considered as
    submanifold rather than a link of two components. 
    Without loss of generality we may assume that
    ${\mathcal F}=(v_1,v_2)$ is smooth and orthogonal.
    Let $\Hat{W}\cong M\times [0,1]$ be the trace of a shift from
    $M$ to $M^{sh}$ along $v_1$. Cut out a small $\varepsilon$-neighbourhood
    $U$ of $Z\times\{\frac{1}{2}\}\subset \Hat{W}$, and define
    $W=\Hat{W}\setminus U$, where $Z=s^{-1}(-1)$. See Fig.1. 

\begin{figure}[htb]
\setlength{\unitlength}{1bp}
\begin{picture}(255,255)(-40,0)
\epsfxsize=10cm
\put(-60,125){\epsfbox{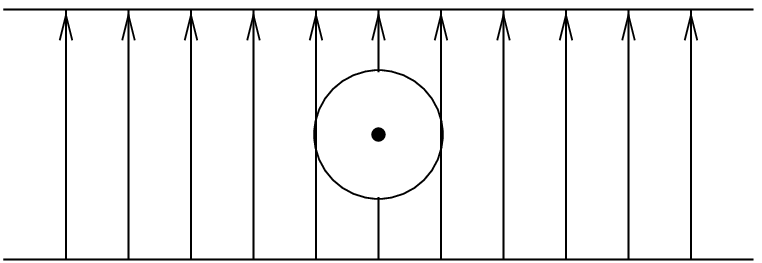}}
\put(-55,170){$\Hat W$}
\put(20,170){$W$}
\put(80,160){$Z$}
\put(85,180){$\partial U$}
\put(177,195){$v_1$}
\end{picture}
\vspace{-4cm}
\caption{}
\end{figure}

       Provide $W\subset\Bbb R^m\times\{0\}\subset\Bbb R^{m+1}$
    with the framing ${\mathcal G}=(e_{m+1},v_2)$, where $e_{m+1}$ is 
    the last vector in the usual base of $\Bbb R^{m+1}$. Let
    $Z\times [-\varepsilon,\varepsilon]\subset M$ be a 
    $\varepsilon$-neighbourhood
    of $Z$. Up to homotopy we may assume
    $s$ maps $M\setminus Z\times [-\varepsilon,\varepsilon]$ to the
    base point $1\in S^1$. We get now a well defined map 
    $\hat{s}:W\longrightarrow S^1$, given by
\[\hat{s}(x,t)=\left\{\begin{array}{r@{\quad:\quad}l}
                                      s(x) & t\ge\frac{1}{2},\\
                                      1 & t\le\frac{1}{2}.
\end{array}\right.\]
    So we obtain a new framing $\hat{s}{\mathcal G}$ of
    $W\subset \Bbb R^{m+1}$. In addition, it holds
    $\partial W=M\sqcup M^{sh}\sqcup\partial U$, where
    $\partial U\cong Z\times S^1$ is the boundary of $U$. We construct now 
    a diffeotopy of $\Bbb R^{m+1}$ which deforms $(W,\hat{s}{\mathcal G})$
    to a framed bordism.

       Let $\nu(M^{sh})$ be the normal vector bundle of
    $M^{sh}\subset\Bbb R^{m+1}$, framed
    by $(v_1,e_{m+1},v_2)$. A homotopy of $\nu(M^{sh})$
           $$F_1:\nu(M^{sh})\times [0,1]\longrightarrow \nu(M^{sh})$$
    is given by rotating $e_{m+1}$
    to $v_1$, $v_1$ to $-e_{m+1}$ and meanwhile keeping $v_2$ fixed.
    Define 
       $$ F_1':M^{sh}\times [0,1]\longrightarrow \Bbb R^{m+1} $$
    by $F_1'(x,t)=(x,-t)$. Let $U_1$ be a $\delta$-neighbourhood of
    $M^{sh}\subset\Bbb R^{m+1}$ with $\delta\ll\epsilon$.
    From $F_1$, $F_1'$ we get an isotopy
    $H_1:U_1\times [0,1]\longrightarrow \Bbb R^{m+1}$, given by
\begin{eqnarray*}
        &  & H_1(x+r_1v_1(x)+r_2e_{m+1}+r_3v_2(x),t) \\ \nonumber
        & =& F_1'(x,t)+r_1F_1(v_1(x),t)+r_2F_1(e_{m+1},t)+
                                  r_3F_1(v_2(x),t), \nonumber
\end{eqnarray*}
   where $x\in M^{sh}$ and 
   $x+r_1v_1(x)+r_2e_{m+1}+r_3v_2(x)\in U_1$. It holds clearly
     $$H_1(M^{sh},1)=M^{sh}\times\{-1\}\subset\Bbb R^m\times\{-1\}.$$
   Since $(e_{m+1},v_2)$ is deformed to $(v_1,v_2)$ and 
   $\hat s|_{M^{sh}}=s$ ($s$ is defined on $M^{sh}$ by identifying
   $M^{sh}$ with $M$ in the natural way),
   $\hat{s}{\mathcal G}|_{M^{sh}}$ is homotoped to $s{\mathcal F}$.

\begin{figure}[htb]
\setlength{\unitlength}{1bp}
\begin{picture}(255,255)(-40,0)
\epsfxsize=10cm
\put(-60,40){\epsfbox{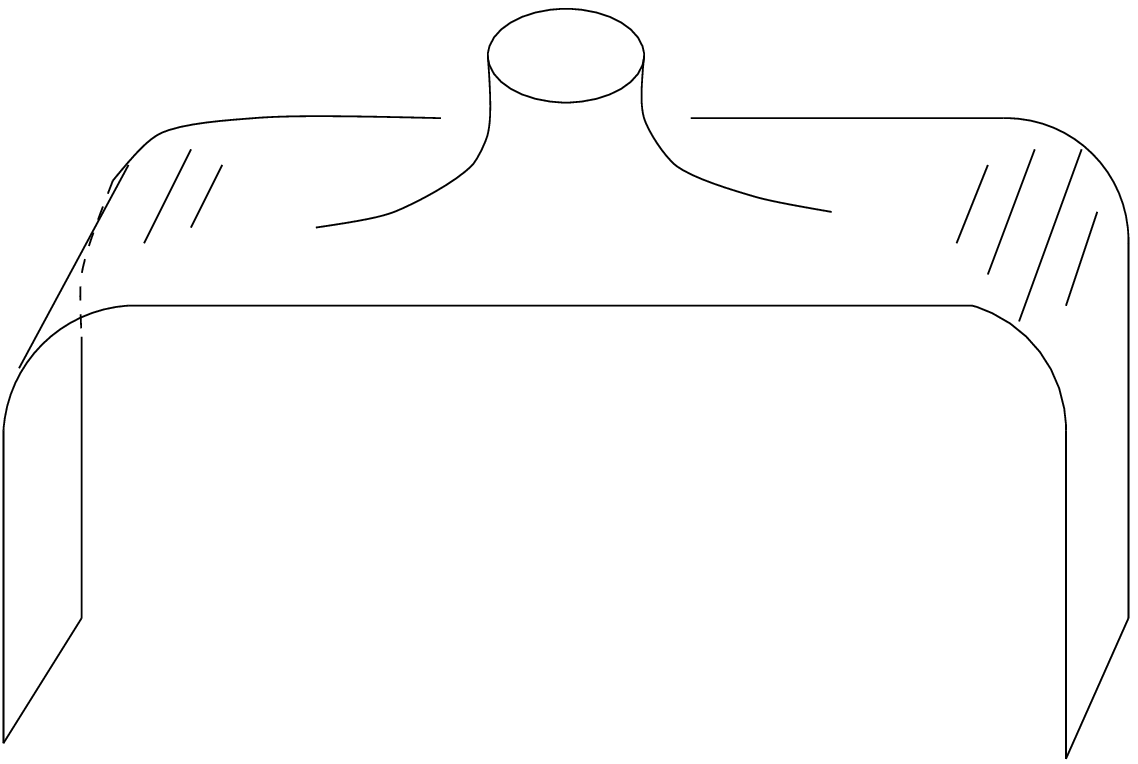}}
\put(74,215){$\partial U$}
\put(-30,60){$M^{sh}$}
\put(180,60){$M$}
\end{picture}
\vspace{-1.5cm}
\caption{}
\end{figure}

        Let $U_0$ be a $\delta$-neighbourhood of $M\subset\Bbb R^{m+1}$.
    Similarly we have an isotopy
        $$ H_0:U_0\times [0,1]\longrightarrow\Bbb R^{m+1} $$ 
    which deforms $M$ to
    $M\times\{-1\}\subset\Bbb R^m\times\{-1\}$ and
    $(e_{m+1},v_2)$ to $(-v_1,v_2)=u{\mathcal F}$. We have used the
    homotopy $F_0:\nu(M)\times [0,1]\longrightarrow \nu(M)$ which rotates
    $e_{m+1}$ to $-v_1$ and keeps $v_2$ fixed.

      Let $\nu(\partial U)$ be the normal vector bundle of
   $\partial U\subset\Bbb R^{m+1}$, framed by $(u_1,e_{m+1},v_2)$,
   where $u_1$ is the normal vector field of $\partial U\subset W$
   pointing inwards. Let
      $$F_2:\nu(\partial U)\times [0,1]\longrightarrow \nu(\partial U)$$
   be the homotopy given by rotating $e_{m+1}$ to $u_1$ and keeping
   $v_2$ fixed; and define
       $$F_2':\partial U\times [0,1]\longrightarrow\Bbb R^{m+1}$$ 
   by $F_2'(x,t)=(x,t)$. From $F_2$ and $F_2'$ we obtain an isotopy
   $H_2:U_2\times [0,1]\longrightarrow \Bbb R^{m+1}$,
   where $U_2$ is a $\delta$-neighbourhood of $\partial U\subset\Bbb R^{m+1}$.
   $H_2$ isotopes $\partial U$ to
   $\partial U\times\{1\}\subset\Bbb R^{m}\times\{1\}$ and
   homotopes $(e_{m+1},v_2)$ to $(u_1,v_2)$. 
   So $\hat{s}{\mathcal G}|_{\partial U}$
   is homotoped to $\hat{s}(u_1,v_2)$. It is not difficult to see that
   $\hat{s}(u_1,v_2)$ is homotopic to the $2\pi$-rotation of the
   $S$-framing $(u_1,v_2)$ of $\partial U\subset\Bbb R^m\times\{1\}$,
   and therefore 
      $$(\partial U,\hat{s}(u_1,v_2))=(Z\times S^1,\hat{s}(u_1,v_2))$$
   is just the fibre-wise embedding of the framed circle $S^1\subset\Bbb R^3$
   representing $\eta:S^3\longrightarrow S^2$.
   It follows $[\partial U,\hat{s}(u_1,v_2)]=\eta_*[Z,{\mathcal F},s]$.

       Now $U_{\delta}=U_0\cup U_1\cup U_2$ is a
   $\delta$-neighbourhood of $\partial W\subset\Bbb R^{m+1}$.
   The isotopies $H_0$, $H_1$ and $H_2$ together define an isotopy
   $H:U_{\delta}\times [0,1]\longrightarrow\Bbb R^{m+1}$.
   According to a well known theorem in Differential Topology $H$ determine 
   a diffeotopy $\widetilde{H}$ of $\Bbb R^{m+1}$ which deforms
   $(W,\hat{s}{\mathcal G})$ to a framed bordism from
   $[(M\sqcup M^{sh}),(u{\mathcal F}\sqcup s{\mathcal F})]$ to
   $\eta_*[Z,{\mathcal F},s]$, see Fig.2.
\end{proof}
\begin{Corollary} (i)\  {\bf [Turaev, 1985]:} 
   It holds $[M,s{\mathcal F}]= \eta_*[Z,{\mathcal F},s]$,
   if ${\mathcal F}$ is the $S$-framing of the submanifold $M$;
   
   (ii)\  $E[M,s{\mathcal F}]-E[M,{\mathcal F}]=E\eta_*[Z,{\mathcal F},s]$,
   where $E$ is the suspension homomorphism.
\end{Corollary}

\begin{proof} (i)\  Let $F$ be a Seifert surface of $M$ giving the
   $S$-framing ${\mathcal F}$. Because $M^{sh}$ is a small shift of $M$
   along ${\mathcal F}$, we have $M^{sh}\pitchfork F=\phi$. Using
   $F$ and $M^{sh}\times [0,1]$ we obtain easily a framed bordism
   from $[(M\sqcup M^{sh}),(u{\mathcal F}\sqcup s{\mathcal F})]$ to
   $[M,s{\mathcal F}]$. The assertion follows now from Proposition 2.1.
   See also Turaev [15].

   (ii)\  We have clearly
\begin{eqnarray}
              E[(M\sqcup M^{sh}),(u{\mathcal F}\sqcup s{\mathcal F})]
              &=& E[M,u{\mathcal F}]+E[M,s{\mathcal F}]  \nonumber \\
              &=& -E[M,{\mathcal F}]+E[M,s{\mathcal F}]. \nonumber
\end{eqnarray}
    The statement follows from Proposition 2.1.
\end{proof}

        In general it is a subtle problem to measure the difference
  $[M,s{\mathcal F}]-[M,{\mathcal F}]$. 
  Hilton-Hopf invariant up to order $3$ are involved.
  Consider the framed link 
  $(M,u{\mathcal F})\sqcup (M^{sh},s{\mathcal F})$ 
  representing an element $\alpha$ in
               $$\pi_m(S^2_1\vee S^2_2)\cong
                   \oplus_{\gamma}\pi_m(S^{q(\gamma)+1}),  $$
  where $\gamma$ runs through a system $\Gamma$ of basic Whitehead products in
  $\iota_1<\iota_2$, and 
  $q(\gamma)$ is the height of $\gamma\in\Gamma$,
  see Hilton [2]. So we have the splitting
  $\alpha=\oplus_{\gamma}\alpha_{\gamma}$ with
  $\alpha_{\gamma}\in \pi_m(S^{q(\gamma)+1})  $.
  Let $\alpha_1$, $\alpha_2$, $\alpha_3$ be the Hilton coefficients of
  $\alpha$ corresponding to 
  $[\iota_1,\iota_2]$, $[\iota_1,[\iota_1,\iota_2]]$ and
  $[\iota_2,[\iota_1,\iota_2]]$ respectively. 
  If we map the wedge $S^2_1\vee S^2_2$ canonically to the sphere $S^2$, 
  then $\iota_1$ and $\iota_2$ are both identified with the identity $\iota$
  of $S^2$, and $\alpha$ is mapped to the element
      $$[(M\sqcup M^{sh}),(u{\mathcal F}\sqcup s{\mathcal F})]
        \in\pi_m(S^2).$$
  The basic Whitehead products in $\iota_1<\iota_2$ are mapped to Whitehead
  products in $\iota$. Because the Whitehead products in $\iota$ with
  weight $>3$ are zero homotopic, we have
\begin{eqnarray}
          &{} & [(M\sqcup M^{sh}),(u{\mathcal F}
               \sqcup s{\mathcal F})] \nonumber \\
          &=&  [M,u{\mathcal F}]+[M^{sh},s{\mathcal F}]
               +[\iota,\iota]_*\alpha_1
          +[\iota,[\iota,\iota]]_*(\alpha_2+\alpha_3). \nonumber
\end{eqnarray}

         Let $\alpha_1'$ be the first nontrivial Hilton-Hopf invariant of
   $[M,{\mathcal F}]$, namely the one corresponding to $[\iota_1,\iota_2]$,
   then it holds (see Hilton [2])
               $$[M,u{\mathcal F}]=-[M,{\mathcal F}]
                 +[\iota,\iota]_*\alpha_1'.$$
   By this and Proposition 2.1 we get 
\begin{Corollary}  It holds
\begin{eqnarray}
                   [M,s{\mathcal F}]-[M,{\mathcal F}]
                   &=& \eta_*[Z,{\mathcal F},s]
                       -[\iota,\iota]_*(\alpha_1+\alpha_1')
                   -[\iota,[\iota,\iota]]_*(\alpha_2+\alpha_3). \nonumber
\end{eqnarray}
   Note that $3[\iota,[\iota,\iota]]=0$ by Jacobi-identity. 
\end{Corollary}

      Consider now a framed link
   $(M_1,{\mathcal F}_1)\sqcup (M_2,{\mathcal F}_2)\subset \mathbb{R}^m$ 
   of codimensions $k_1, k_2\ge 2$.
   The bordism group $FL_m^{k_1,k_2}$ of such links is isomorphic to
   $\pi_m (S^{k_1}\vee S^{k_2})$. We try to understand the role of the framings
   ${\mathcal F}_1$, ${\mathcal F}_2$ in the Hilton splitting of
   $(M_1,{\mathcal F}_1)\sqcup (M_2,{\mathcal F}_2)$.

      Let $i=1,2$. A Seifert surface of $M_i$ is a compact oriented submanifold
   $F_i\subset\Bbb R^m$ with boundary $M_i$. If $M_i$ has a framed
   Seifert surface $F_i$ such that the framing of $M_i$ as the boundary of
   $F_i$ is homotopic to the original framing ${\mathcal F}_i$, 
   then we say ${\mathcal F}_i$
   is an $S$-framing of $M_i$ and $M_i$ is $S$-framed.
   We call $F_i$ a suitably framed Seifert surface.
   Note that two $S$-framings must not be homotopic.

\begin{Proposition} Let
   $(M_1,{\mathcal F}_1)\sqcup (M_2,{\mathcal F}_2)\subset \mathbb{R}^m$
    be an $S$-framed link representing
   $\alpha\in\pi_m (S^{k_1}\vee S^{k_2})$, and let $F_1$, $F_2$ be
   the corresponding suitably framed Seifert surfaces.

     (i)\  Up to involution it holds
          $[M_1\pitchfork F_2]=[F_1\pitchfork M_2].$
          An involution is an isomorphism $u$ of the target group
          with $u\circ u=id$.

     (ii)\  Let $\alpha=\oplus_{\gamma}\alpha_{\gamma}$ be the
              Hilton splitting of $\alpha$. Up to involution 
              $\alpha_{[\iota_1,\iota_2]}$ is given by  
              $[M_1\pitchfork F_2]$, all other Hilton coefficients 
              $\alpha_{\gamma}$ are zero.
\end{Proposition}
\begin{proof} The desired framed bordism in the first part is given by
     $F_1\pitchfork F_2$.

          Let $Z=M_1\pitchfork F_2$ and let $U\cong Z\times D^{k_1}$
      be an open tubular neighbourhood of $Z\subset F_2$.
      Frame $\partial U$ so that $F_2\setminus U$ 
      gives rise to a framed bordism between $M_2$ and $\partial U$.
      Because $F_2\setminus U$ is disjoint from $M_1$,
      $M_1\sqcup M_2$ is framed bordant to $M_1\sqcup\partial U$, 
      see Fig.3.

\begin{figure}[htb]
\setlength{\unitlength}{1bp}
\begin{picture}(255,255)(-40,0)
\epsfxsize=11cm
\put(-60,140){\epsfbox{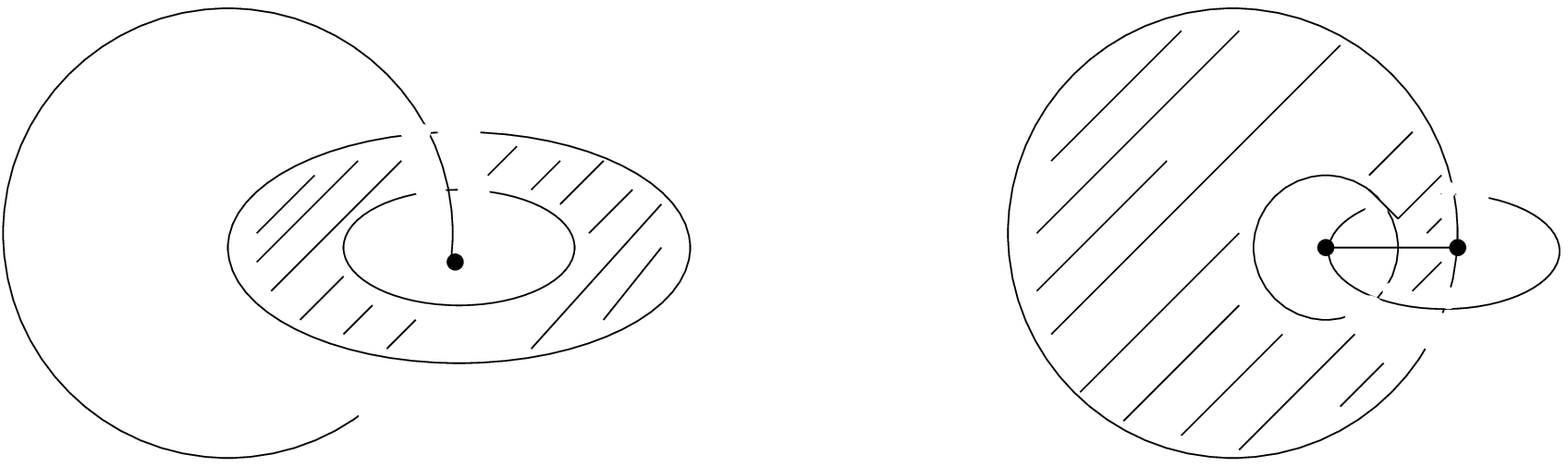}}
\put(-50,195){$M_1$}
\put(30,165){$F_2$}
\put(50,210){$M_2$}
\put(17,182){$Z$}
\put(38,182){$\partial U$}
\put(177,210){$F_1$}
\put(128,210){$M_1$}
\put(195,205){$\partial U'$}
\put(192,185){$Z'$}
\put(235,180){$Z$}
\put(240,200){$\partial U$}
\end{picture}
\vspace{-4.5cm}
\caption{}
\end{figure}

        If $U$ is small enough, then it holds
     $F_1\pitchfork \bar U\cong Z\times [0,1]$ with
     $Z=Z\times\{0\}$. Let $v$ be the normal vector field of
     $Z\subset F_1\pitchfork \bar U$ pointing inwards.
     It's easy to see that $Z'=Z\times\{1\}=F_1\pitchfork \partial U$
     is just a small shift of $Z$ along $v$. Let
     $U'$ be a small tubular neighbourhood of
     $Z'\subset F_1$ such that $U'\pitchfork \partial U=Z'$.
     Frame $\partial U'$ so that $F_1\setminus U'$ gives rise to a 
     framed bordism between $M_1$ and $\partial U'$. 
     Because $F_1\setminus U'$ is disjoint from $\partial U$ it follows
     $[M_1\sqcup \partial U]=[\partial U'\sqcup\partial U]$,
     see Fig.3 again.

          In addition, $\partial U'\sqcup\partial U
      =Z\times S^{k_2-1}\sqcup Z\times S^{k_1-1}$ 
      is just the fibre-wise embedding
      of the standard framed Hopf link
      $ S^{k_2-1}\sqcup S^{k_1-1}\subset\mathbb{R}^{k_1+k_2-1}$
      (at least up to involution of the framings) into a small
      tubular neighbourhood of the framed intersection 
      $Z\subset \mathbb{R}^m$. It follows
               $$[M_1\sqcup M_2]=[\iota_1,\iota_2]_*[Z]$$
      at least up to involution. Since the Hilton splitting of 
      $[\iota_1,\iota_2]_*[Z]$ has the form
             $$ 0+0+[Z]+0+\cdots, $$ 
      the assertion follows.
\end{proof}

\begin{Corollary}  Let $\gamma$ be a basic Whitehead product in
     $\iota_1<\iota_2$ of weight $\ge 3$, and let
     $M_1\sqcup M_2\subset\mathbb{R}^{q(\gamma)+1}$
     be a framed link representing $\gamma$. At least one component
     of this link is not $S$-framed.
\end{Corollary}
\begin{proof} If this is not the case then according to the above result
     we have $H_{\gamma}[M_1\sqcup M_2]=0$, a contradiction to the fact
     $H_{\gamma}[M_1\sqcup M_2]=\pm 1\in \pi_{q(\gamma)+1}(S^{q(\gamma)+1})
      \cong Z$, where $H_{\gamma}$ is the Hilton homomorphism
     corresponding to the basic Whitehead product $\gamma$.
\end{proof}

\section{isomorphisms between $FL^2_{m,r}$ and $L^2_{m,r+1}$}

         Sanderson [12] gave an isomorphism
    $\theta:FL^2_{m,r}\longrightarrow L^2_{m,r+1}$. We construct here for each
    subset $\sigma\subset\{1,2,\cdots,r\}$ such an 
    isomorphism $\theta_{\sigma}$.

          Let $L=(M_1,{\mathcal F}_1)\sqcup\cdots\sqcup (M_r,{\mathcal F}_r)
     \subset\mathbb{R}^m $ be a framed link of codimension $2$.
     For $i\in\sigma$ define $\Hat M_i=-M_i^{sh}$, where
     $-M_i^{sh}$ is the negative oriented $M_i^{sh}$. 
     Note that a framed submanifold
     in $\Bbb R^m$ is canonically oriented. For $i\not\in\sigma$
     consider the framed intersection
     $Z_i=M_i\pitchfork \widetilde{F}_i$, where $\widetilde{F}_i$
     is a Seifert surface of $M_i^{sh}$.
     Let $U_i$ be a small tubular neighbourhood of 
     $Z_i\subset \widetilde{F}_i$.
     If we orient $\partial U_i$ as the boundary of $U_i$ then
     $\widetilde{F}_i\setminus U_i$ is a Seifert surface of
     $(-M_i^{sh})\sqcup (-\partial U_i)$. For $i\not\in\sigma$ we define
     $\Hat M_i=\partial U_i$ and $M_{r+1}=\sqcup_{i=1}^r \Hat M_i$.
     $\theta_{\sigma}$ is given by the assignment
                 $$\theta_{\sigma}[L]=
                     [M_1\sqcup\cdots\sqcup M_r\sqcup M_{r+1}]_{or}.$$
\begin{Lemma} $\theta_{\sigma}:FL^2_{m,r}\longrightarrow L^2_{m,r+1}$
      is a well defined homomorphism.
\end{Lemma}
\begin{proof} To prove this let
     $W_1\sqcup\cdots\sqcup W_r\subset\mathbb{R}^m\times [0,1]$
     be a framed bordism between
       $$ L=\sqcup_{i=1}^r (M_i,{\mathcal F}_i),\ \ \ \ \  
          L'=\sqcup_{i=1}^r (M'_i,{\mathcal F}'_i).    $$ 
     Take $\widetilde{W}_i=W_i^{sh}$
     and let $\widetilde{F}_i$, $\widetilde{F}'_i$ be Seifert surfaces
     of $M_i^{sh}$ and ${M'_i}^{sh}$. We define for $i\not\in\sigma$
                    $$\Hat W_i=\widetilde{W}_i\cup \widetilde{F}_i\times \{0\}
                         \cup \widetilde{F}'_i\times \{1\}$$
     and orient $\Hat W_i$ so that its orientation is in agreement with
     $\widetilde {W_i}$.
     Let $F_{\Hat W_i}$ be a Seifert surface of $\Hat W_i$,
     consider $\Hat Z_i=W_i\pitchfork F_{\Hat W_i}$ with boundary
     $\partial \Hat Z_i=Z_i\times\{0\}\sqcup Z'_i\times\{1\}$.
     Cut out a small tubular neighbourhood $\Hat U_i$ of
     $\Hat Z_i\subset F_{\Hat W_i}$ and orient $\partial \Hat U_i$ so
     that $\partial \Hat U_i$ is an oriented bordism from $\partial U_i$
     to $\partial U'_i$. Define
             $$W_{r+1}=(\sqcup_{i\not\in\sigma} \partial \Hat U_i)
                  \sqcup (-\sqcup_{i\in\sigma}W_i^{sh}). $$
     $W_1\sqcup\cdots\sqcup W_r\sqcup W_{r+1}$ is clearly an oriented
     bordism between $\theta_{\sigma}(L)$ and
     $\theta_{\sigma}(L')$. It is clear that $\theta_{\sigma}$ 
     respects the addition.
\end{proof}

            To prove that $\theta_{\sigma}$ is an isomorphism we construct now
     a homomorphism $\zeta_{\sigma}:L^2_{m,r+1}\longrightarrow FL^2_{m,r}$
     and show it's in fact the inverse of $\theta_{\sigma}$.

          Let $\bar{\sigma}=\sigma\cup\{r+1\}$, and let
    $L=M_1\sqcup\cdots\sqcup M_r\sqcup M_{r+1}\subset\mathbb{R}^m$
    be an oriented link of codimension $2$. Take a Seifert surface
    $F_{\bar{\sigma}}$ of $M_{\bar{\sigma}}=\sqcup_{i\in\bar{\sigma}} M_i$.
    Denote by $u_1$ the normal vector field of
    $M_{\bar{\sigma}}\subset F_{\bar{\sigma}}$ pointing outwards and by
    $u_2$ the normal vector field of $F_{\bar{\sigma}}$ determined by the
    orientation. 
    For $i\in\sigma$ take ${\mathcal F}_i=(u_1,u_2)|_{M_i}$ as a framing
    of $M_i$.  If $i\not\in\sigma$ consider the intersection
    $Z_i=M_i\pitchfork F_{\bar{\sigma}}$ and let
    $Z_i\times [-1,1]\subset M_i$ be a small tubular neighbourhood of
    $Z_i\subset M_i$ such that the 
    positive direction of $[-1,1]$ is in agreement
    with $u_2|_{Z_i}$. We assume here that $M_i$ intersects $F_{\bar{\sigma}}$
    perpendicularly. Let $(v_1^S,v_2^S)$ be the $S$-framing of $M_i$
    determined by the orientation. 
    Define ${\mathcal F}_i=(v_1,v_2)$
    to be the $2\pi$-rotation of $(v_1^S,v_2^S)$ in
    $Z_i\times [-1,1]$, see (1) and (2) in \S 2.
    By doing this we have provided $M_i$ with the framing ${\mathcal F}_i$
    for $1\le i\le r$.
    Now $\zeta_{\sigma}:L_{m,r+1}^2\longrightarrow FL_{m,r}^2$ 
    is defined by the following assignment
         $$\zeta_{\sigma}[M_1\sqcup\cdots\sqcup M_r\sqcup M_{r+1}]_{or}=
           [(M_1,{\mathcal F}_1)\sqcup\cdots\sqcup (M_r,{\mathcal F}_r)].$$
\begin{Lemma} $\zeta_{\sigma}$ is a well defined homomorphism.
\end{Lemma}
\begin{proof} Let $W_1\sqcup\cdots\sqcup W_r\sqcup W_{r+1}\subset
     \mathbb{R}^m\times [0,1]$ be an oriented bordism between
        $$ L=M_1\sqcup\cdots\sqcup M_r\sqcup M_{r+1},
           \hskip1cm L'=M_1'\sqcup\cdots\sqcup M_r'\sqcup M'_{r+1} $$ 
     and let $F'_{\bar{\sigma}}$ be a Seifert surface of
     $M'_{\bar{\sigma}}=\sqcup_{i\in\bar{\sigma}} M'_i$. Consider
            $$W=W_{\bar{\sigma}}\cup F_{\bar{\sigma}}\times \{0\}
                       \cup F'_{\bar{\sigma}}\times \{1\}, $$
     where $W_{\bar{\sigma}}$ is the disjoint union of the $W_i$'s with
     $i\in \bar{\sigma}$. We orient $W$ so that its 
     orientation coincides with the one of $W_{\bar{\sigma}}$.   
     Smooth $W$ and let $F_W$ be an oriented
     Seifert surface of $W$. $F_W$ induces canonically 
     a framing ${\mathcal G}_i$
     on $W_i$ for each $i\in\sigma$. If $i\not\in\sigma$ consider
             $$\Hat W_i=W_i\cup F_i\times \{0\}
                  \cup F'_i\times \{1\},$$
     where $F_i$, $F_i'$ are Seifert surfaces of $M_i$ and $M_i'$,
     and again we orient $\Hat W_i$ so that its orientation is in
     agreement with the one of $W_i$.
     Let ${\mathcal G}^S_i$ be the restriction of the $S$-framing of $\Hat W_i$
     on $W_i$. Consider now the intersection $\Hat Z_i=W_i\pitchfork F_W$
     with boundary
                 $$\partial \Hat Z_i=Z_i\times\{0\}\sqcup Z_i'\times\{1\},$$
     and denote by $\Hat Z_i\times [-1,1]\subset W_i$ a small tubular
     neighbourhood such that the positive direction of $[-1,1]$ coincides with
     the normal vector field of $F_W\subset\Bbb R^{m+1}$ determined by the
     orientation. For $i\not\in\sigma$ let ${\mathcal G}_i$ be the framing
     of $W_i$ given by the $2\pi$-rotation of ${\mathcal G}^S_i$
     in $\Hat Z_i\times [-1,1]$. Clearly
              $$(W_1,{\mathcal G}_1)\sqcup\cdots\sqcup (W_r,{\mathcal G}_r)$$
      is a framed bordism from $\zeta_{\sigma}[L]$ to $\zeta_{\sigma}[L']$,
      and $\zeta_{\sigma}$ evidently respects the addition.
\end{proof}
\begin{Theorem} $\theta_{\sigma}:FL_{m,r}^2\longrightarrow L_{m,r+1}^2$
      is an isomorphism, in fact its inverse is $\zeta_{\sigma}$.
\end{Theorem}
\begin{proof}  (1)\ \  $\zeta_{\sigma}\circ\theta_{\sigma}=id$. Consider
            \begin{eqnarray}
              \zeta_{\sigma}\circ\theta_{\sigma}
              [(M_1,{\mathcal F}_1)\sqcup\cdots\sqcup (M_r,{\mathcal F}_r)]&=&
             \zeta_{\sigma}
             [M_1\sqcup\cdots\sqcup M_r\sqcup M_{r+1}]_{or} \nonumber  \\
             &=& [(M_1,{\mathcal F}'_1)\sqcup\cdots
             \sqcup (M_r,{\mathcal F}'_r)].  \nonumber
          \end{eqnarray}
    For $M_{r+1}$ and the framings ${\mathcal F}_i'$ see the definitions of
    $\theta_{\sigma}$ and $\zeta_{\sigma}$. 
    Denote by $M_i\times [0,\varepsilon]$
    the trace of a small $\varepsilon$-shift from $M_i$ to $M_i^{sh}$. Orient
             $$F_{\bar{\sigma}}=(\sqcup_{i\not\in\sigma} U_i)\sqcup
                  (\sqcup_{i\in\sigma} M_i\times [0,\varepsilon])$$
     so that $F_{\bar{\sigma}}$ is a Seifert surface of
     $M_{\bar{\sigma}}$. According to the definition of $\zeta_{\sigma}$
     we see easily that up to homotopy it holds 
     ${\mathcal F}'_i={\mathcal F}_i$ for
     $i\in\sigma$. In addition, for $i\not\in\sigma$ the intersection
     $M_i\pitchfork F_{\bar{\sigma}}=M_i\pitchfork U_i$ 
     is exactly the submanifold
     $Z_i=M_i\pitchfork \widetilde{F}_i$, because $U_i$ is a small tubular
     neighbourhood of $Z_i\subset\widetilde{F}_i$, see the definition of
     $\theta_{\sigma}$. This implies that ${\mathcal F}_i$ and 
     ${\mathcal F}'_i$ are
     homotopic, since both ${\mathcal F}_i$ and ${\mathcal F}'_i$ 
     are essentially the
     $2\pi$-rotation of the $S$-framing of $M_i$ 
     in a small tubular neighbourhood
     of $Z_i\subset M_i$. It follows 
     $\zeta_{\sigma}\circ\theta_{\sigma}=id$.

       (2)\ \ $\theta_{\sigma}\circ\zeta_{\sigma}=id$. Consider
               \begin{eqnarray}
               \theta_{\sigma}\circ\zeta_{\sigma}
               [M_1\sqcup\cdots\sqcup M_r\sqcup M_{r+1}]_{or}&=&
               \theta_{\sigma}
               [(M_1,{\mathcal F}_1)\sqcup\cdots\sqcup (M_r,{\mathcal F}_r)]
                               \nonumber  \\
           &=& [M_1\sqcup\cdots\sqcup M_r\sqcup M'_{r+1}]_{or}.
                               \nonumber
\end{eqnarray}
      For $M'_{r+1}$ and the framings see the definitions of $\theta_{\sigma}$
      and $\zeta_{\sigma}$. Let $F_{\bar{\sigma}}$ be an oriented Seifert
      surface of $M_{\bar{\sigma}}$. Cut out a small tubular neighbourhood
      $V_{\sigma}$ of $M_{\sigma}\subset F_{\bar{\sigma}}$ to get
      $F'_{\bar{\sigma}}$ with
         $$\partial F'_{\bar{\sigma}}=\partial (F_{\bar{\sigma}}\setminus
           V_{\sigma})=M_{r+1}\sqcup (\sqcup_{i\in\sigma} M_i^{sh}).$$
      For $i\not\in\sigma$ let $\widetilde{F}_i$ be an oriented Seifert
      surface of $M_i^{sh}$. According to the definitions of $\theta_{\sigma}$
      and $\zeta_{\sigma}$ we may assume
                   $$Z_i=M_i\pitchfork \widetilde{F}_i
                             =M_i\pitchfork F'_{\bar{\sigma}}$$
      for $i\not\in\sigma$, because $\partial\widetilde{F}_i=M_i^{sh}$
      is a shift of $M_i$ along the framing 
      ${\mathcal F}_i$ ( ${\mathcal F}_i$ is given
      by the $2\pi$-rotation of the 
      $S$-framing in a small tubular neighbourhood
      of $M_i\pitchfork  F'_{\bar{\sigma}}\subset M_i$), and because
      $M_i\pitchfork\widetilde{F}_i$ is just where the
      $2\pi$-rotation of the $S$-framing takes place. 
      In addition, we may assume
      $M_i$ intersects $\widetilde{F}_i$ perpendicularly. This implies
      $U_i\subset F'_{\bar{\sigma}}$ after a small isotopy of
      $F'_{\bar{\sigma}}$ fixing boundary, 
      where $U_i\subset \widetilde{F}_i$
      is a small tubular neighbourhood of $Z_i=M_i\pitchfork \widetilde{F}_i$.
      Define
               $$\Hat F_{\bar{\sigma}}=F'_{\bar{\sigma}}
                    \setminus (\cup_{i\not\in\sigma} U_i)$$
      which is a Seifert surface of
                 $$M_{r+1}\sqcup(\sqcup_{i\in\sigma} M_i^{sh})
                      \sqcup(\sqcup_{i\not\in\sigma} -\partial U_i)
                      =M_{r+1}\sqcup (-M_{r+1}'). $$
      Because $\Hat F_{\bar{\sigma}}$ is disjoint from all $M_i$, 
      $1\le i\le r$,
      there is an embedding $W_{r+1}\subset\mathbb{R}^m\times[0,1]$
      of $\Hat F_{\bar{\sigma}}$ such that $W_{r+1}$ is an oriented bordism
      from $M_{r+1}$ to $M'_{r+1}$ and such that $W_{r+1}$ is disjoint from
      $M_i\times [0,1]\subset\Bbb R^m\times [0,1]$ for all $1\le i\le r$.
      So we obtain an oriented  bordism
               $$M_1\times [0,1]\sqcup\cdots\sqcup M_r\times [0,1]
                    \sqcup W_{r+1}\subset \mathbb{R}^m\times [0,1]$$
      from $M_1\sqcup\cdots\sqcup M_r\sqcup M_{r+1}$ to
      $M_1\sqcup\cdots\sqcup M_r\sqcup M'_{r+1}$.
      It follows $\theta_{\sigma}\circ\zeta_{\sigma}=id$.
\end{proof}
\begin{Example} If $\sigma$, $\sigma'\subset\{1,\cdots,r\}$ are different
    then $\theta_{\sigma}\not=\theta_{\sigma'}$ in general.
    Without loss of generality we assume $1\in\sigma$ and $1\not\in\sigma'$.
    Consider the framed Hopf link
             $$L=S^1_1\sqcup S^1_2\sqcup\phi\sqcup
                 \cdots\sqcup\phi\subset \mathbb{R}^3.$$
     Let $L_{\sigma}$, $L_{\sigma'}$ be the oriented links representing
     $\theta_{\sigma}[L]$ and $\theta_{\sigma'}[L]$ respectively.
     It is not difficult to verify that the linking number between the second
     and the last components of $L_{\sigma}$ is $\pm 1$, and the linking
     number between the second and the last components of $L_{\sigma'}$
     is $0$. This shows $\theta_{\sigma}\not=\theta_{\sigma'}$.
     In the case $\sigma=\{1,\cdots,r\}$ it holds $\theta_{\sigma}=\theta$,
     where $\theta$ is the isomorphism given by Sanderson [12].
     Some computations may be simplified by suitably choosing
     $\theta_{\sigma}$ or $\zeta_{\sigma}$.
     In next section we will choose $\zeta_{\sigma}=\zeta_{\phi}$,
     namely $\sigma=\phi$.
\end{Example}

\section{computation of the forgetful homomorphism}

      Let $f:FL_{m,2}^2\longrightarrow L_{m,2}^2$ be the 
    forgetful homomorphism
    given by forgetting the framings. To compute $f$ we compute
    the composition
       $$f'=\zeta_{\sigma}\circ f:FL_{m,2}^{2}\longrightarrow
            L_{m,2}^{2}\stackrel{\cong}\longrightarrow FL_{m,1}^{2}$$
    where and throughout this section $\sigma=\phi$.

       Let $M^{m-2}\subset\Bbb R^m$ be a framed or an oriented submanifold.
   Denote by $M^S$ the same submanifold but provided with the $S$-framing.
   Let $M_1\sqcup M_2\subset\Bbb R^m$ be a framed link representing
   an element in $FL_{m,2}^2$. According to the definition of $\zeta_{\sigma}$
   and from Corollary 2.2 we see
              $$f'[M_1\sqcup M_2]_{fr}=\eta_*[M_1^S\pitchfork F_2]_{fr}$$
   at least up to involution, where $\eta:S^3\longrightarrow S^2$ 
   is the Hopf map and $F_2$ is a Seifert surface of $M_2$.
   So we need only to compute $M^S_1\pitchfork F_2$.

      It is easily seen that $f'\circ\gamma_*=0$ for $\gamma=\iota_1$
   and $\gamma=\iota_2$. So let $\gamma$ be a basic Whitehead product 
   in $\iota_1<\iota_2$ of weight $\ge 2$,
   and $\eta_{\gamma}:S^{q(\gamma)+1}\longrightarrow S^3$ be the map
   determined by $M^S_1(\gamma)\pitchfork F_2(\gamma)$, where
   $M_1(\gamma)\sqcup M_2(\gamma)$ is a framed link representing
   $\gamma$ and $F_2(\gamma)$ is a framed Seifert surface of
   $M_2(\gamma)$. Given $[Z]\in\pi_m(S^{q(\gamma)+1})$, 
   by fibre-wise embedding it follows easily that
   $f'\circ\gamma_*[Z]=\eta_{\gamma}[Z]$ which is represented
   by $Z\times (M^S_1(\gamma)\pitchfork F_2(\gamma))$. Therefore
   we only need to compute $M^S_1(\gamma)\pitchfork F_2(\gamma)$.
  
       For $\gamma=[\iota_1,\iota_2]$ we see easily that
    $M^S_1(\gamma)\pitchfork F_2(\gamma)$ is a framed point,
    this means $\eta_{\gamma}$ is the identity up to sign in this case.  
    Let $\gamma_w=[\iota_2,[\iota_2\cdots[\iota_2,[\iota_1,\iota_2]]\cdots]]$
    be a basic Whitehead product of weight $w\ge 3$, clearly $q(\gamma_w)=w$.
    The following framed link in $\mathbb{R}^{w+1}$ represents $\gamma_w$
\begin{eqnarray}
             M_1(w)&=&S^1_w\times S^1_{w-1}\times\cdots
                             \times S^1_3\times S^1_2, \nonumber  \\
             M_2(w)&=&S^1_w\times S^1_{w-1}\times\cdots
                             \times S^1_3\times S^1_1\sqcup \nonumber  \\
                             & & S^1_w\times S^1_{w-1}\times\cdots
                                 \times S^2\sqcup \nonumber  \\
                             & & \cdots\cdots\cdots\sqcup \nonumber  \\
                             & & S^1_w\times S^{w-2}\sqcup \nonumber  \\
                             & & S^{w-1} \nonumber  \\
                             &=& N_{2,1}\sqcup\cdots\sqcup N_{2,w-1}, \nonumber
\end{eqnarray}
     where $S^1_2\sqcup S^1_1$, $S^1_i\sqcup S^{i-1}$ are
     $S$-framed Hopf links, $3\le i\le w$, and all products are given by
     fibre-wise embeddings. $M_1(w)$ is clearly $S$-framed, 
     but $M_2(w)$ not, according to Corollary 2.5. In addition, we can assume
     $M_1(w)\subset \mathbb{R}^w\times\{0\}\subset\mathbb{R}^{w+1} $.

          Consider now the oriented submanifolds of $M_1(w)$
                    $$Z_i=S^1_w\times\cdots\times S^1_{i+1}\times\{pt\}
                         \times S^1_{i-1}\cdots\times S^1_3\times S^1_2,$$
    $2\le i\le w$ and $\{pt\}$ is a set of a single point. 
    $\cup_{i=2}^w Z_i\subset M_1(w)$ is a transversally
    immersed submanifold. Using the classical trick in Fig.4 we can 
    successively dissolve
    the multi-points to get an embedded submanifold $Z(w)\subset M_1(w)$
    of codimension $1$. Let $v_1$, $v_2$ be the normal vector fields of
    $Z(w)\subset M_1(w)$ and of $M_1(w)\subset\Bbb R^w\times\{0\}$
    respectively. Define ${\mathcal F}_w=(v_1,v_2|_{Z(w)})$.

\begin{figure}[htb]
\setlength{\unitlength}{1bp}
\begin{picture}(255,255)(-40,0)
\epsfxsize=12cm
\put(-85,150){\epsfbox{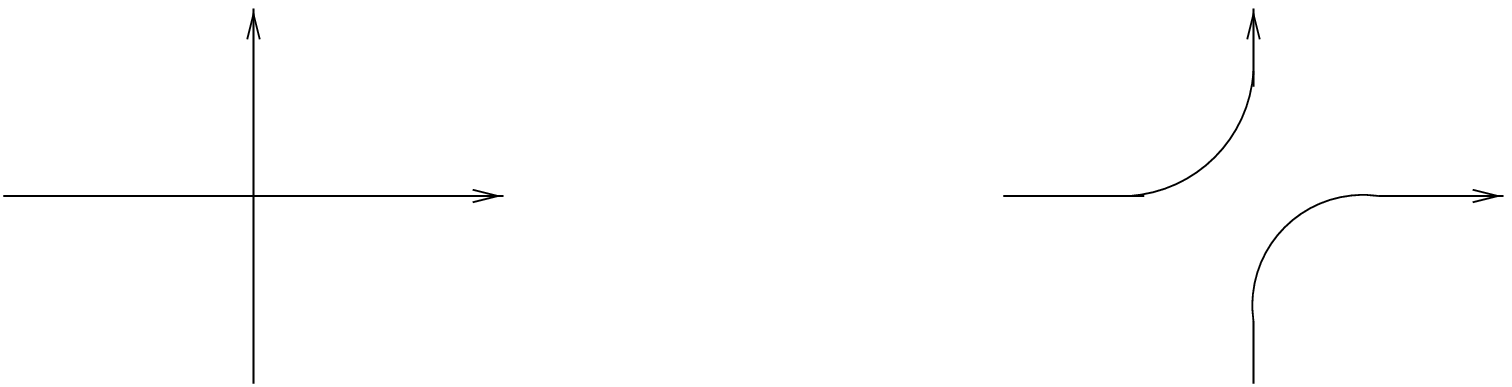}}
\end{picture}
\vspace{-4.7cm}
\caption{}
\end{figure}

\begin{figure}[htb]
\setlength{\unitlength}{1bp}
\begin{picture}(255,255)(-40,0)
\epsfxsize=8cm
\put(-25,-10){\epsfbox{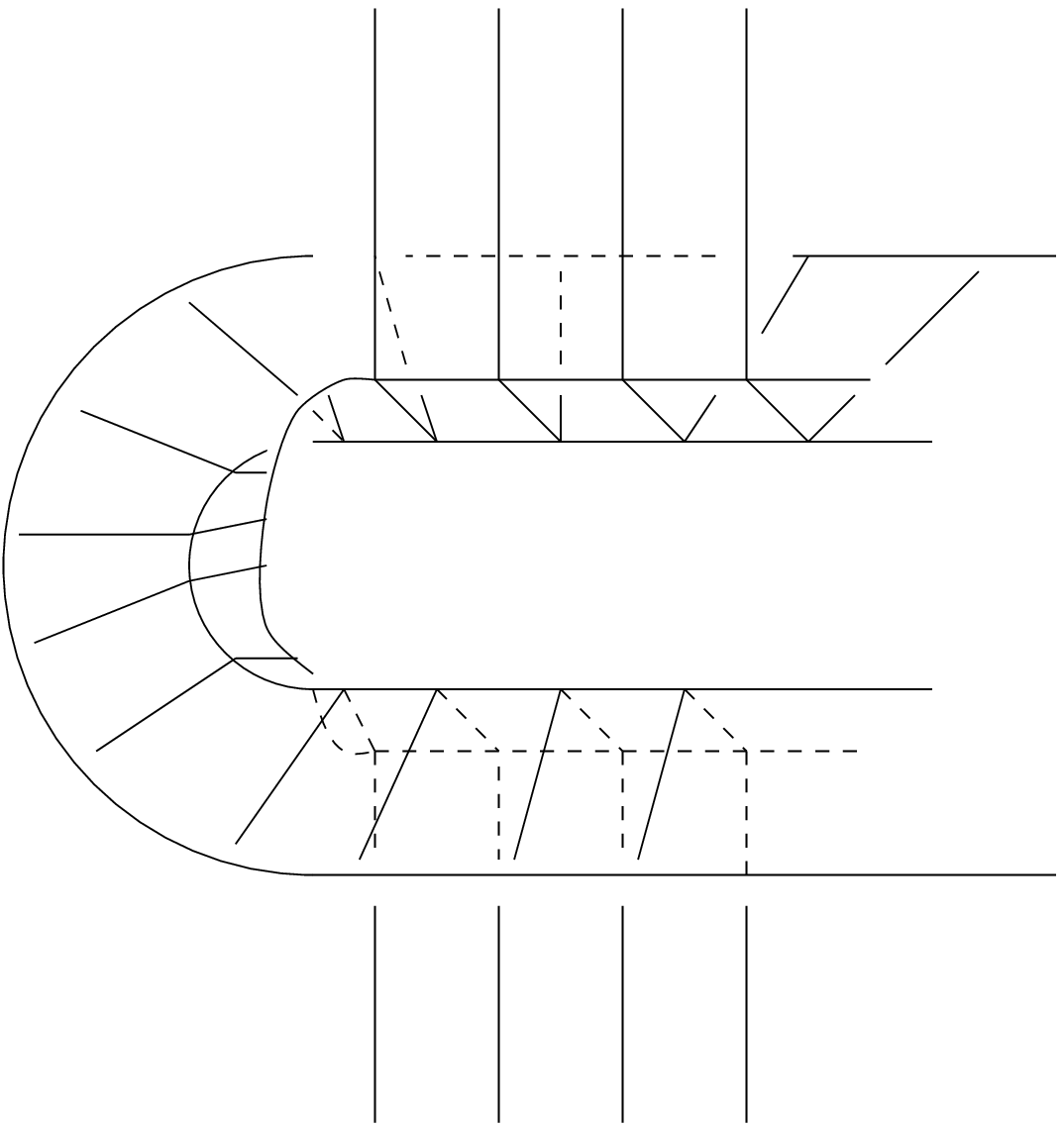}}
\put(-15,169){$N_1^{sh}$}
\put(10,150){$F_1^*$}
\put(85,200){$F_2^*$}
\put(30,200){$N_2^{sh}$}
\end{picture}
\vspace{0.5cm}
\caption{}
\end{figure}

       For the proof of Proposition 4.1 we need to consider the following
     situation. Let 
             $$ N^{m-2}=N_1\sqcup N_2\subset\Bbb R^m $$ 
     be an oriented closed submanifold and let $F_1$, $F_2$ be 
     Seifert surfaces of $N_1$ 
     and $N_2$ respectively. We construct now a Seifert surface
     $F$ of $N$ from $F_1$ and $F_2$. Cut out tubular neighborhoods of
     $N_1\subset F_1$ and $N_2\subset F_2$ to get $F_1^*$ and $F_2^*$,
     and cut out a tubular neighbourhood $U(C)$ of 
         $$ C=F_1^*\cap F_2^*\subset F_1^*\cup F_2^*. $$
     Then we can use the trick in Fig.5 to sew $F_1^*\setminus U(C)$ and
     $F_2^*\setminus U(C)$ together to get $F^*$ with
     $\partial F^*=N_1^{sh}\sqcup N_2^{sh}$, where $N_1^{sh}$ is essentially
     a small shift of $N_1$ along the framing given by the $2\pi$-rotation
     of the $S$-framing of  $N_1$ in a tubular neighbourhood of
     $N_1\pitchfork F_2\subset N_1$, similarly $N_2^{sh}$.
     Let $T_1$, $T_2$ be the traces of these shifts.
     $F^*\cup T_1\cup T_2$ is oriented, we can make it smooth and get a
     Seifert surface $F$ of $N$. Note that, shown in Fig.5
     is the locus near $F_1\pitchfork N_2$, for the locus
     $N_1\pitchfork F_2$ it is completely similar; away from
     $F_1\pitchfork N_2$ and $N_1\pitchfork F_2$ the method shown
     in Fig.4 applies. For details see [16], p.14--17.

\begin{Proposition} (a)\ Let $\gamma_w$, $M_1(w)$, $M_2(w)$, $Z(w)$ 
     and ${\mathcal F}_w$ be as above, and assume $w\ge 3$.
     There exists a Seifert surface $F_2(w)$
     of $M_2(w)$ such that the following holds

        (i)\  $M_1(w)\pitchfork F_2(w)=Z(w)$,

        (ii)\  $[M^S_1(w)\pitchfork F_2(w)]=
                         \pm E[Z(w),{\mathcal F}_w].$

    (b)\  Up to sign the framed submanifold
    $M^S_1(w)\pitchfork F_2(w)\subset\Bbb R^{w+1}$
    represents the map
         $$\eta_{\gamma}=(E\eta)\circ (E^2\eta)\circ\cdots\circ (E^{w-3}\eta)
              \circ (E^{w-2}\eta):S^{w+1}\longrightarrow S^3.$$
\end{Proposition}
\begin{proof}
     (a)\  The case $w=3$ is easy. Assume inductively that
     the assertion is true for $\gamma_{w-1}$, and let $F_2(w-1)$ be a
     Seifert surface of $M_2(w-1)$ with the desired property.

         By fibre-wise embedding we get a Seifert surface
    $S^1_w\times F_2(w-1)$ of $N_{2,1}\sqcup\cdots\sqcup N_{2,w-2}$.
    $N_{2,w-1}=S^{w-1}$ bounds a ball $D^w$. From $D^w$ and
    $S^1_w\times F_2(w-1)$ we obtain a Seifert surface $F_2(w)$ of
    $M_2(w)$. To compute $M_1(w)\pitchfork F_2(w)$ we look at
             \begin{eqnarray}
                 Z' &=& M_1(w)\pitchfork S^1_w\times F_2(w-1) \nonumber \\
                    &=& S^1_w\times (M_1(w-1)\pitchfork F_2(w-1)) \nonumber \\
                    &=& S^1_w\times Z(w-1), \nonumber \\
                Z'' &=& M_1(w)\pitchfork D^w \nonumber \\
                    &=& \{pt\}\times S^1_{w-1}\times\cdots
                        \times S^1_3\times S^1_2
                        \nonumber \\
                    &=& Z_w. \nonumber
\end{eqnarray}
    Considered in $M_1(w)$ we have the transversal intersection
                  $$Q=Z'\pitchfork Z''=\{pt\}\times Z(w-1) $$
    of codimension $2$. Because $Q$ is disjoint from the boundaries of
    $S^1_w\times F_2(w-1)$ and $D^w$ we see that in the construction
    of $F_2(w)$ we have just dissolved $Q$ as Fig.4. This means
    $M_1(w)\pitchfork F_2(w)=Z(w)$. Part (i) follows. 
    Because $M_1(w)$ is $S$-framed, at least up to sign
    ${\mathcal F}_w$ is the framing of 
    $M^S_1(w)\pitchfork F_2(w)=Z(w)$,
    part (ii) follows. 

    (b)\  Up to involution we have
                 $$[M^S_1(w)\pitchfork F_2(w)]=
                      [F_1(w)\pitchfork M^S_2(w)],$$
     see Proposition 2.4. Consider
\begin{eqnarray*}
     Z'(w) &=& F_1(w)\pitchfork M_2(w) \nonumber \\
           &=& F_1(w)\pitchfork N_{2,1}  \nonumber \\
           &=& S^1_w\times S^1_{w-1}\times\cdots\times S^1_3\times D_2
               \pitchfork \\ \nonumber
           & &   S^1_w\times S^1_{w-1}\times
               \cdots\times S^1_3\times S^1_1, \nonumber \\
           &=& S^1_w\times S^1_{w-1}\times\cdots
               \times S^1_3\times \{pt\}, \nonumber
\end{eqnarray*}
     where $D_2$ is a disk with boundary $S^1_2$. Let $(v_1,v_2)$ be the
     $S$-framing of $M_2(w)$ and $v_3$ be the normal vector field of
     $F_1(w)\subset\Bbb R^{w+1}$. 
     Define ${\mathcal F}'_w=(v_1,v_2,v_3)|_{Z'(w)}$.
     So we have
             $$[Z'(w),{\mathcal F}'_w]=[F_1(w)\pitchfork M^S_2(w)],$$
     In addition, $(v_1,v_2)|_{N_{2,1}}$ is given by the 
     $2\pi$-rotations of the $S$-framing
     $(u_1,u_2)$ of $N_{2,1}$ in the tubular neighborhoods of all
             $$S^1_w\times\cdots\times S^1_{i+1}\times\{pt\}\times S^1_{i-1}
                  \times\cdots\times S^1_3\times S^1_1\subset N_{2,1},$$
     $3\le i\le w$ and $\{pt\}$ denotes a set of a single point. 
     We will get $(Z'(w),{\mathcal F}'_w)$ when we take a regular value
     in $S^3$ of the map
       $$(E\eta)\circ (E^2\eta)\circ\cdots\circ (E^{w-3}\eta)
                \circ (E^{w-2}\eta):S^w\longrightarrow S^3$$
     and perform the Pontryagin-Thom construction. The statement follows.
\end{proof}     
\begin{Proposition} Let $\gamma=[\iota_2,[\iota_2,\cdots[\iota_2,
     [\iota_1,[\iota_1,\cdots[\iota_1,\iota_2]\cdots]]]\cdots]]$
     be a basic Whitehead product in $\iota_1<\iota_2$ of weight $w\ge 3$,
     and let $M_1\sqcup M_2\subset\Bbb R^{w+1}$ be a framed link
     representing $\gamma$. The following map
        $$ \eta_{\gamma}=(E\eta)\circ (E^2\eta)\circ\cdots\circ (E^{w-3}\eta)
                      \circ (E^{w-2}\eta):S^{w+1}\longrightarrow S^3$$
     is represented by $M^S_1\pitchfork F_2$ up to sign,
     where $F_2$ is a Seifert surface of $M_2$.
\end{Proposition}
\begin{proof} Assume $\iota_i$ appears $w_i$-times in $\gamma$,
     $i=1,2$. The following framed link $M_1\sqcup M_2$
     in $\Bbb R^{w+1}$ representing $\gamma$
\begin{eqnarray*}
     M_1&=&S^1_w\times S^1_{w-1}\times\cdots \times S^1_{w_1+2}
            \times S^1_{w_1+1}\cdots
            \times S^1_3\times S^1_2\sqcup \nonumber  \\
        & & S^1_w\times S^1_{w-1}\times\cdots
            \times S^1_{w_1+2}
            \times S^1_{w_1+1}\cdots \times S^2\sqcup  \nonumber  \\
        & & \cdots\cdots\cdots\sqcup  \nonumber  \\
        & & S^1_w\times S^1_{w-1}\times \cdots
            \times S^1_{w_1+2}\times S^{w_1}
                        \nonumber  \\
        &=& N_{1,1}\sqcup\cdots\sqcup N_{1,w_1} \nonumber \\
\end{eqnarray*}
\begin{eqnarray*} 
   M_2 &=& S^1_w\times S^1_{w-1}\times\cdots\times S^1_{w_1+3}
             \times S^1_{w_1+2}\times\cdots
             \times S^1_3\times S^1_1\sqcup  \nonumber  \\
        & & S^1_w\times S^1_{w-1}\times\cdots\times
             S^1_{w_1+3}\times S^{w_1+1}\sqcup  \nonumber  \\
        & & \cdots\cdots\cdots \sqcup  \nonumber  \\
        & & S^1_w\times S^{w-2}\sqcup  \nonumber  \\
        & & S^{w-1} \nonumber  \\
        &=& N_{2,1}\sqcup\cdots\sqcup N_{2,w_2}. \nonumber
\end{eqnarray*}
    where $S_2^1\sqcup S^1_1$, $S_3^1\sqcup S^2,\cdots$,
    $S_w^1\sqcup S^{w-1}$ are usually framed Hopf links and all products
    are given by fibre-wise embeddings.

        According to Proposition 4.1 we can assume
\begin{eqnarray}
                Z &=&  F_1\pitchfork M_2=F_1\pitchfork N_{2,1} \nonumber \\
                  &=&  S^1_w\times S^1_{w-1}\times \cdots
                              \times S^1_{w_1+2}\times Z(w_1+1), \nonumber
\end{eqnarray}
    where $Z(w_1+1)$ is given by dissolving the multi-points of the
    following immersion iteratedly
\begin{eqnarray*}
      &  \cup_{i=1,i\not=2}^{w_1+1} S^1_{w_1+1}\times\cdots\times S^1_{i+1}
           \times\{pt\}\times S^1_{i-1}\times\cdots\times S^1_3\times S^1_1
                           \\           
      &  \subset S^1_{w_1+1}\times\cdots\times S^1_3\times S^1_1.
\end{eqnarray*}
    Let $Z'(w_1+1)=S^1_{w_1+1}\times\cdots\times S^1_3\times\{pt\}$,
    framed by $(v_1,v_2,v_3)$ as in the proof of Proposition 4.1, part (b). 
    So we have a framed
    bordism $(W',{\mathcal G}')$ from $Z(w_1+1)$ to $Z'(w_1+1)$ (at least up
    to involution) with ${\mathcal G}'=(\bar{v}_1,\bar{v}_2,\bar{v}_3)$. 
    Consider the fibre-wise embedding
          $$W=S^1_w\times\cdots\times S^1_{w_1+2}\times W'\subset
                    \Bbb R^{w+1}\times [0,1]. $$
    We obtain a framing ${\mathcal G}$ of $W$ by performing the
    $2\pi$-rotations of $(\bar{v}_2,\bar{v}_3)$ in tubular neighborhoods
    of all
           $$S^1_w\times\cdots\times S^1_{i+1}\times\{*\}\times S^1_{i-1}
                  \times\cdots\times S^1_{w_1+2}\times W'\subset W,$$
    $w_1+2\le i\le w$. $(W,{\mathcal G})$ is a framed bordism from
    $F_1\pitchfork M^S_2$ to $Z'(w)$,
    according to the proof of Proposition 4.1, part (b), 
    $Z'(w)$ represents
            $$ (E\eta)\circ (E^2\eta)\circ\cdots\circ (E^{w-3}\eta)
                \circ (E^{w-2}\eta):S^{w+1}\longrightarrow S^3.$$
    The assertion follows by Lemma 2.4.
\end{proof}

       Note that in case $m=4$ it holds $f'\circ\gamma_*\not=0$
    for $\gamma=[\iota_1,\iota_2]$, $[\iota_1,[\iota_1,\iota_2]]$
    and $[\iota_2,[\iota_1,\iota_2]]$, and $f'\circ\gamma_*=0$
    for $\gamma=\iota_1$ and $\iota_2$. Because $\zeta_{\phi}$
    is an isomorphism, the same holds if we replace $f'$ by $f$.
    So we recover the corresponding results of Sanderson [12].
    
       Let $\Gamma$ be a system of basic Whitehead products 
    in $\iota_1<\iota_2$, and
    $\gamma=[\alpha,\beta]\in\Gamma$ be such that the weights of
    $\alpha$ and $\beta$  are greater than $1$. Take $\gamma'$ to be one of
    $\alpha$, $\beta$ and $\gamma$, and let
    $M_1(\gamma')\sqcup M_2(\gamma')\subset \Bbb R^{q(\gamma')+1}$
    be a framed link representing $\gamma'$. Denote by
    $\eta_{\gamma'}:S^{q(\gamma')+1}\longrightarrow S^3$ the map
    given by $[M^S_1(\gamma')\pitchfork F_2(\gamma')]_{fr}$,
    where $F_2(\gamma')$ is a Seifert surface of $M_2(\gamma')$.

\begin{Proposition} Under the above notations and up to sign we have
           $$ \eta_{\gamma}=[\eta_{\alpha},\eta_{\beta}]:
                 S^{q(\alpha)+q(\beta)+1}=
                S^{q(\gamma)+1}\longrightarrow S^3.$$
\end{Proposition}
\begin{proof} Surely, we may take
    $M_1(\gamma)\sqcup M_2(\gamma)\subset \Bbb R^{q(\gamma)+1}$
    to be the following
\begin{eqnarray}
     M_1(\gamma)&=&S^{q(\alpha)}\times M_1(\beta)\sqcup
                   S^{q(\beta)}\times M_1(\alpha),\nonumber \\
     M_2(\gamma)&=&S^{q(\alpha)}\times M_2(\beta)\sqcup
                   S^{q(\beta)}\times M_2(\alpha),\nonumber
\end{eqnarray}
    where $S^{q(\alpha)}\sqcup S^{q(\beta)}\subset\Bbb R^{q(\gamma)+1}$
    is the usually framed Hopf link and all products are given by fibre-wise
    embeddings. Let $F_i(\alpha)$, $F_i(\beta)$ be Seifert surfaces
    of $M_i(\alpha)$ and $M_i(\beta)$ respectively, $i=1,2$, then we get a
    Seifert surface
            $$F_i(\gamma)=S^{q(\alpha)}\times F_i(\beta)\sqcup
                                        S^{q(\beta)}\times F_i(\alpha) $$
    of $M_i(\gamma)$ by fibre-wise embeddings. This implies
\begin{eqnarray}
    M^S_1(\gamma)&=&S^{q(\alpha)}\times M^S_1(\beta)\sqcup
                    S^{q(\beta)}\times M^S_1(\alpha), \nonumber \\
    M^S_1(\gamma)\pitchfork F_2(\gamma)&=&
                    S^{q(\alpha)}\times
                    (M^S_1(\beta)\pitchfork  F_2(\beta)) \sqcup
                    S^{q(\beta)}\times
                   (M^S_1(\alpha)\pitchfork  F_2(\alpha)).  \nonumber
\end{eqnarray}
    Because by induction $M^S_1(\alpha)\pitchfork  F_2(\alpha)$,
    $M^S_1(\beta)\pitchfork F_2(\beta)$ represent the maps
    $\eta_{\alpha}$ and $\eta_{\beta}$ respectively, the assertion
    $\eta_{\gamma}=[\eta_{\alpha},\eta_{\beta}]$ follows.
\end{proof}

       Combining Propositions 4.1, 4.2 and 4.3 we get a complete
    computation of $f'$ and therefore the forgetful homomorphism $f$.
    Given $\alpha=\oplus_{\gamma}\alpha_{\gamma}\in\pi_m(S^2\vee S^2)$,
    $f'(\alpha)$ is the sum of $f'\circ\gamma_*(\alpha_{\gamma})$,
    and all the homomorphisms $f'\circ\gamma_*$ 
    are computed above. To get the decomposition
    $\alpha=\oplus_{\gamma}\alpha_{\gamma}$ we may use the method 
    in the author's paper [17].
        We end this paper with a question which U. Kaiser mentioned to me.
    Given a framed link 
    $(M_1,{\mathcal F}_1)\sqcup (M_2,{\mathcal F}_2)\subset\Bbb R^m$
    of codimension $2$, denote by ${\mathcal F}^S_i$ the
    $S$-framing of $M_i$, $i=1,2$. There is a map
    $s_i:M_i\longrightarrow S^1$ such that
    ${\mathcal F}_i\backsimeq s_i{\mathcal F}^S_i$.
    Let $q_i:S^1\longrightarrow S^1$ be a map of degree $q_i$ and
    define $s'_i=q_i\circ s_i$. The assignment
       $$[(M_1,{\mathcal F}_1)\sqcup (M_2,{\mathcal F}_2)]\longmapsto
         [(M_1,s'_1{\mathcal F}^S_1)\sqcup (M_2,s'_2{\mathcal F}^S_2)] $$
    determines a well defined homomorphism
    $(q_1,q_2)_*:FL_{m,2}^2\longrightarrow FL_{m,2}^2$.
    The question is
\begin{quote}
          How can one describe the homomorphism $(q_1,q_2)_*$
     by using the Hilton splitting of $FL_{m,2}^2\cong\pi_m(S^2\vee S^2)$ ?
\end{quote}

\vskip .3cm

{\small
{\it
Address

Department of Mathematics

Bar-Ilan University

52900 Ramat Gan, Israel}

e-mail: wangji@macs.biu.ac.il  
}
\end{document}